\documentclass{article}
\usepackage{amsmath,amsfonts,amsthm,amssymb}
\usepackage{graphicx}

\newtheorem{thm}{Theorem}[section]
\newtheorem{lem}[thm]{Lemma}
\theoremstyle{remark}
\newtheorem{rmk}{Remark}
\newtheorem{eg}{Example}
\newcommand{\B}{\mathcal B}
\newcommand{\M}{\mathcal M}

\begin{document}

\title{Combinatorial proofs of inverse relations and \\ log-concavity for Bessel numbers}
\author{%
Hyuk Han\thanks{Department of Mathematics Education, Seoul National University, Seoul 151-748, South Korea. Email:\,hyuk0916@snu.ac.kr.} 
\quad{and}\quad
Seunghyun Seo\thanks{Department of Mathematics, Brandeis University, Waltham, MA 02454, USA. Email:\,shseo@brandeis.edu. Research supported by the Post-doctoral Fellowship Program of Korea Science and Engineering Foundation (KOSEF).}}
\maketitle
\begin{abstract}
Let the Bessel number of the second kind $B(n,k)$ be the number of  set partitions of $[n]$ into $k$ blocks of size one or two, and let the Bessel number of the first kind $b(n,k)$ be the coefficient of $x^{n-k}$ in~$-y_{n-1}(-x)\,$, where $y_{n}(x)$ is the $n$th Bessel polynomial. In this paper, we show that Bessel numbers satisfy two properties of Stirling numbers: The two kinds of  Bessel numbers are related by inverse formulas, and both Bessel numbers of the first kind and the second kind form log-concave sequences. By constructing sign-reversing involutions, we prove the inverse formulas. We review Krattenthaler's injection for the log-concavity of Bessel numbers of the second kind, and give a new explicit injection for the log-concavity of signless Bessel numbers of the first kind.
\end{abstract}

%%%%%%%%%%%
%%%%%%%%%%%
\section{Introduction}
\subsection{Main Results}
For nonnegative integer $n$ and $k$\,, let $\B(n,k)$ be the set of partitions of $[n]:=\{1,\ldots , n\}$ into $k$ blocks of size one or two. Let  $B(n,k)$ denote the number of elements in $\B(n,k)\,$. The Bessel polynomials \cite{Bo, B, CS, KF} are the polynomial solutions $y_{n}(x)$ of the second-order differential equations
\begin{equation*}
x^{2}\,y_{n}''+(2x+2)\,y_{n}'=n(n+1)\,y_n \,, \quad y_{n}(0)=1\,.
\end{equation*}
Let $a(n,k)$ denote the coefficient of $x^{n-k}$ in $y_{n-1}(x)\,$. In fact, $a(n,k)$ equals the cardinality of the set $\B(2n-k-1,n-1)\,$. Set $b(n,k):=(-1)^{n-k}\,a(n,k)\,$.

In this paper we deduce following two identities.
\begin{equation}  \label{form:bessel-inv}
\sum_{k=0}^{n}\,B(n,k)\, b(k,l)= \delta_{n,l}
\quad \mbox{and}
\quad \sum_{k=0}^{n}\,b(n,k)\, B(k,l) = \delta_{n,l}\,.
\end{equation}
These identities are reminiscent for the inverse formulas of Stirling numbers \cite[Prop 1.4.1]{St}. By constructing sign-reversing involutions, we give combinatorial proofs for both identities.

Second, we show the log-concavity of $\{\,a(n,k)\,\}_{1\leq k \leq n}\,$, i.e., that
\begin{equation*}
a(n,k-1)\cdot a(n,k+1) \leq a(n,k)^2 \,,
\end{equation*}
which is an analogue of the log-concavity of the signless Stirling numbers of the first kind.  We prove this result by constructing an explicit injection from the set $\B(2n-k,n-1)\times \B(2n-k-2,n-1)$ to the set $\B(2n-k-1,n-1)\times \B(2n-k-1,n-1)\,$, which provides a combinatorial proof of the log-concavity of the sequence.

%%%%%%%%%%%%%%%%%%%%
\subsection{Notations and Backgrounds}
We review some standard terminology, which will be used throughout the paper.

\subsubsection*{Stirling numbers}
For a nonnegative integer $n\,$, let $[n]$ denote the set $\{1,2,\ldots,n\}$\,, and let $[0]=\emptyset$\,. A {\em partition} of $[n]$ is a collection of nonempty pairwise disjoint subsets whose union is $[n]\,$. These disjoint subsets of $[n]$ are called {\em blocks} of the set partition. The number of partitions of $[n]$ with $k$ blocks $B_1 , \ldots, B_k$ is called the {\em Stirling number of the second kind} $S(n,k)\,$. By convention, we put $S(0,0)=1$\,. We denote by $S_n$ the set of all permutations of $[n]\,$. Let $c(n,k)$ be the number of permutations in $S_n$ that have $k$ cycles. Define $s(n,k)$ to be $(-1)^{n-k}\,c(n,k)\,$. This number $s(n,k)$ is called the {\em Stirling number of the first kind}. It is well-known that Stirling numbers are related to one another by the following {\em inverse formulas} \cite[Prop 1.4.1]{St}:
\begin{equation} \label{form:st-inv}
\sum_{k=0}^{n}\,S(n,k)\,s(k,l) = \delta_{n,l}
\quad \mbox{and}
\quad \sum_{k=0}^{n}\,s(n,k)\,S(k,l) = \delta_{n,l}\,.
\end{equation}
Moreover, both Stirling numbers satisfy the {\em log-concavity} property, i.e., that
\begin{equation} \label{form:st-log}
c(n,k-1)\cdot c(n,k+1) \leq c(n,k)^2
\quad \mbox{and}
\quad S(n,k-1)\cdot S(n,k+1) \leq S(n,k)^2 \, .
\end{equation}
There are many proofs for (\ref{form:st-inv}) and (\ref{form:st-log}) with both algebraic and combinatorial ways. See \cite{St, SW,W}.

\subsubsection*{Bessel numbers}
Recall that the Bessel polynomials are the polynomial solutions $y_{n}(x)$ of the second-order differential equations
\begin{equation} \label{equ:bessel_diff}
x^{2}\,y_{n}''+(2x+2)\,y_{n}'=n(n+1)\,y_n
\end{equation}
satisfying the initial condition $y_{n}(0)=1\,$.
Bochner \cite{Bo} seems first to have realized that these polynomials are closely related to the Bessel functions. Krall and Frink \cite{KF} have considered  the system of Bessel polynomials in connection with certain solutions of the wave equations. Moreover, Burchnall \cite{B} has developed more properties of the Bessel polynomials in detail. Recently, Choi and Smith \cite{CS} have investigated coefficients of Bessel polynomials from a combinatorial point of view.

From the differential equation (\ref{equ:bessel_diff}) one can  easily derive the formula:
\begin{equation*}
y_{n}(x)=\sum_{k=0}^{n}\,\dfrac{(n+k)!}{2^{k}\,k!\,(n-k)!}\,  x^k \,.
\end{equation*}
Let $a(n,k)$ denote the coefficient of $x^{n-k}$ in $y_{n-1}(x)\,$. Set  $b(n,k):=(-1)^{n-k}\,a(n,k)\,$, i.e.,
\begin{equation} \label{form:1st-bessel}
b(n,k):=\begin{cases}
\,(-1)^{n-k}\,\dfrac{(2n-k-1)!}{2^{n-k}\,(n-k)!\,(k-1)!}\,, & \mbox{if $1 \leq k \leq n\,$,}\\
\,0\,,&\mbox{if $1 \leq n < k\,$.}\end{cases}
\end{equation}
We call the number $b(n,k)$ a {\em Bessel number of the first kind}, and $a(n,k)$ a {\em signless Bessel number of the first kind}. By convention, we put $a(0,k)=b(0,k)=\delta_{0,k}$\,.

For nonnegative integers $n$ and $k$, define $\B(n,k)$ to be the set of partitions of $[n]$ into $k$ nonempty blocks of each size one or two. A block $B$ is called a {\em singleton} if $|B|=1\,$, and a {\em pair} if $|B|=2\,$. Note that if $\pi \in \B(n,k)\,$, then $\pi$ has exactly $2k-n$ singletons and $n-k$ pairs. Set $B(n,k):=|\B(n,k)|\,$. To choose $n-k$ pairs, we should choose $2n-2k$ elements from $[n]$ and pairing the chosen $2n-2k$ elements, so the number $B(n,k)$ is given by
\begin{equation}\label{form:2nd-bessel}
B(n,k) :=\begin{cases}
\,\dfrac{n!}{2^{n-k}\,(n-k)!\,(2k-n)!}\,, & \mbox{if $\lceil {n}/{2}\rceil \leq k \leq n\,$,}\\
\,0\,,&\mbox{otherwise.}\end{cases}
\end{equation}
We call the number $B(n,k)$ a {\em Bessel number of the second kind}.
\begin{rmk}
It is easily checked that $a(n,k)=B(2n-k-1,n-1)\,$. So we can consider $a(n,k)$ as the cardinality of the set $\B(2n-k-1,n-1)\,$.
\end{rmk}

\subsubsection*{Matching number}
Let $G$ be a loopless graph with $n$ vertices. A $k${\em-matching} in a graph $G$ is a set of $k$ edges, no two of which have a vertex in common. A matching $\alpha$ {\em saturates} a vertex $x\,$, or $x$ is said to be {\em saturated} under $\alpha\,$, if an edge in $\alpha$ is incident with $x\,$; otherwise, $x$ is {\em unsaturated} under $\alpha\,$. The complete graph $K_n$ with $n$ vertices is a simple graph with each pair of whose distinct vertices is adjacent. Let $\M(n,k)$ denote the set of $k$-matchings in complete graph $K_n$ and $m(n,k)$ be the cardinality of $\M(n,k)\,$. We call $m(n,k)$ a {\em matching number}. From now on, assume that all vertices of $K_n$ are labeled with the integers in $[n]\,$. Under this assumption, a partition of $[n]$ into $k$ nonempty blocks, each of size at most two, can be identified with an $(n-k)$-matching in $K_n$ canonically. For example, the partition $\pi=\{\,\{1\}, \{2,6\}, \{3\}, \{4,5\}\,\}$ in $\B(6,4)\,$, corresponds to the matching $\alpha=\{\,\{2,6\}, \{4,5\}\,\}$ in $\M(6,2)\,$. Hence $B(n,k)=m(n,n-k)\,$. Whether $\{a,b\}$ is an edge or a pair, we will write it such that $a<b\,$.

It is well known (see \cite{G,K,LP}) that  the matching numbers $\{\,m(n,k)\,\}_{0\leq k \leq n}$ form a log-concave sequence, i.e., that
\begin{equation} \label{form:match-log}
m(n,k-1)\cdot m(n,k+1) \leq m(n,k)^2 \,,
\end{equation}
which implies that  the Bessel numbers of the second kind $B(n,k)$ is also log-concave. In particular, Krattenthaler \cite{K} has given a neat combinatorial proof of (\ref{form:match-log}).
We review Krattenthaler's proof in Section~\ref{subsec:log2}.  

%%%%%%%%%%%%%%%%%%%%%%%%%
%%%%%%%%%%%%%%%%%%%%%%%%%
\section{Inverse formulas for the Bessel numbers}
In this section we deduce the inverse formulas for the Bessel numbers (\ref{form:bessel-inv}) by the generating function technic.

For indeterminates $x$ and $t$, let
\begin{equation*}
\sum_{n=0}^{\infty}\,f_{n}(x)\,\dfrac{t^n}{n!} := \Big(\,1 + t + \dfrac{t^2}{2!}\,\Big)^x \,.
\end{equation*}
Note that $\{\,f_{n}(x)\,\}$ is a sequence of polynomials of binomial type, i.e.,
\begin{equation*}
f_{n}(x+y)=\sum_{k=0}^{n}\,\binom{n}{k}\,f_{k}(x)\,f_{n-k}(y)\,.
\end{equation*}
For details on polynomials of binomial type, see \cite[pp.~167--213]{RM}.

\begin{lem} \label{lem:B}
For any nonnegative integer $n\,$, we have
\begin{equation*}
f_{n}(x) = \sum_{k=0}^{n}\,B(n,k)\,[x]_k \,,
\end{equation*}
where $[x]_{0}=1$ and $[x]_k = x(x-1)\cdots (x-k+1)\,$, for \,$k\geq1$\,.
\end{lem}

\begin{proof} 
We have 
\begin{eqnarray*}
f_{n}(x)
&=&\Big[\dfrac{t^{n}}{n!}\Big]\Big(\,1+t+\dfrac{t^2}{2!}\,\Big)^x \\
&=&\Big[\dfrac{t^{n}}{n!}\Big]\sum_{k\geq 0}\, \binom{x}{k}\,\Big(\,t+\dfrac{t^2}{2!}\,\Big)^k \\
&=&\sum_{k\geq 0}\, \Big[\dfrac{t^{n}}{n!}\Big] \dfrac{(t+\frac{t^2}{2!})^k}{k!} \cdot [x]_{k} \\
&=&\sum_{k\geq 0}\, B(n,k) \cdot [x]_{k}\,.
\end{eqnarray*}
The last equation holds due to the exponential formula \cite[p.~81]{W}\,.
\end{proof}

\begin{lem} \label{lem:b}
For any nonnegative integer $n\,$, we have
\begin{equation*}
[x]_{n} = \sum_{k=0}^{n}\,b(n,k)\,f_{k}(x)\,.
\end{equation*}
\end{lem}

\begin{proof} 
From the formula (\ref{form:2nd-bessel}) and Lemma \ref{lem:B}, we have
\begin{equation*}
f_{n}(x)  = \sum_{k=0}^n\, \frac{n!}{2^{n-k}\,(n-k)!\,(2k-n)!}\,[x]_k \,.
\end{equation*}
Multiplying both sides of this equation by ${2^n}/{n!}$ yields
\begin{eqnarray*}
\frac{2^n}{n!}\,f_{n}(x) 
& = &\sum_{k=0}^n\, \frac{1}{(n-k)!\,(2k-n)!}\,2^k \, [x]_k\\
& = &\sum_{k=0}^n\, \binom{k}{n-k}\, \frac{2^k}{k!}\,[x]_k \,.
\end{eqnarray*}
Wilf \cite[p.~169]{W} shows, using the Lagrange Inversion Formula, that if two sequences $\{a_n\}$ and $\{b_n\}$ are related by
\begin{equation*}
b_n = \sum_{k}\, \binom{k}{n-k}\,a_k \,,
\end{equation*}
then we have the inversion
\begin{equation*}
n\,a_n = \sum_{k}\, \binom{2n-k-1}{n-k}\, (-1)^{n-k}\, k\, b_k \,.
\end{equation*}
So we have
\begin{equation*}
n\,\frac{2^n}{n!}\,[x]_{n}=\sum_{k}\,\binom{2n-k-1}{n-k}\,(-1)^{n-k}\,k\,\frac{2^k}{k!}\,f_{k}(x)\,.
\end{equation*}
Multiplying both sides of this equation by $2^{-n}\,{(n-1)!}$ yields that
\begin{equation*}
[x]_{n}=\sum_{k}\, (-1)^{n-k}\, \frac{(2n-k-1)!}{2^{n-k}\,(n-k)!\,(k-1)!}\,f_{k}(x)\,.
\end{equation*}
From the formula (\ref{form:1st-bessel}), we have
\begin{equation*}
[x]_{n}=\sum_k\, b(n,k)\,f_{k}(x) \,,
\end{equation*}
which completes the proof. 
\end{proof}

Note that both $[x]_n$ and $f_{n}(x)$ are monic polynomials of degree $n$\,. By Lemma \ref{lem:B} and \ref{lem:b}, we have the following {\em inverse formulas} for the Bessel numbers.
\begin{thm}
Let $b(n,k)$ be the Bessel number of the first kind and $B(n,k)$ be the Bessel number of the second kind. Then two Bessel numbers are related to as follows:
\begin{eqnarray} \label{equ:inv}
\sum_{k=0}^{n}\,B(n,k)\,b(k,l) &=& \delta_{n,l}\,, \label{form:1st-inv}\\
\sum_{k=0}^{n}\,b(n,k)\,B(k,l) &=& \delta_{n,l}\,. \label{form:2nd-inv}
\end{eqnarray}
\end{thm}
In Section \ref{sec:combi-inv}, we prove these two inverse formulas (\ref{form:1st-inv}) and (\ref{form:2nd-inv}) combinatorially.

%%%%%%%%%%%%%%%%%%%%%%%%
%%%%%%%%%%%%%%%%%%%%%%%%
\section{Two involutions for the inverse formulas} \label{sec:combi-inv}
Recall that a partition of $[n]$ into $k$ blocks of size one or two can be identified with an $(n-k)$-matching in $K_n$. So we can describe the Bessel numbers in terms of the matching numbers $m(n,k)$ as follows:
\begin{equation*} 
B(n,k)=m(n,n-k) \quad \mbox{and} \quad
b(n,k)=(-1)^{n-k}\,a(n,k)=(-1)^{n-k}\,m(2n-k-1,n-k) \,. 
\end{equation*}

%%%%%%%%%%%%%%%%%%%%%%%%%%
\subsection{An involution for the first inverse formula}
\begin{thm}\label{thm:1ortho}
For all nonnegative integers $n$ and $l$\,, we have
$$\sum_{k=0}^{n}\,B(n,k)\,b(k,l) = \delta_{n,l} \,.$$
\end{thm}
\begin{proof}
Consider the complete graph $K_{2n-l-1}$. For arbitrary $k$\,, let $U$ be the set of ordered pairs $(\alpha,\beta)$ such that
\begin{enumerate}
\item $\alpha$ is an $(n-k)$-matching in $K_{2n-l-1}$ whose saturated vertices are in $[n]$\,; and
\item $\beta$ is a $(k-l)$-matching in $K_{2n-l-1}$ such that $\alpha$ and $\beta$ saturate no vertex in common.
\end{enumerate}
Define the sign of an element $(\alpha,\beta)\in U$ by ${\rm sgn}(\alpha,\beta)=(-1)^{|\beta|}$\,, where $|\beta|$ denotes the number of edges in $\beta$. Then we have
\begin{equation*}
\sum_{(\alpha,\beta)\in U}{\rm sgn}(\alpha,\beta)
=\sum_{k}\,(-1)^{k-l}\,m(n,n-k)\,\,m(2k-l-1,k-l) = \sum_{k}\,B(n,k)\,b(k,l) \,.
\end{equation*}
To prove the theorem, it suffices to find a sign-reversing involution $I_1$ on $U$ having exactly one fixed point with sign $+1$\,.

Consider the union of $\alpha$ and $\beta$\,. Define a linear order on $\alpha \cup \beta$ as follows: For edges $\{a,b\}$ and $\{c,d\}$ in $\alpha \cup \beta$\,, $\{a,b\} < \{c,d\}$ if and only if $b<d$\,. Note that this order is the same as the colex order (see \cite[p.~8]{SW}). With respect to this order, take the smallest edge $e$ of $\alpha \cup \beta$ and move $e$ from $\alpha$ to $\beta\,$, if $e$ is in $\alpha$\,; otherwise, move $e$ from $\beta$ to $\alpha$\,. This defines the needed involution $I_1$ on $U$\,.

Since there are $n-l$ edges in $\alpha \cup \beta$\,, for the smallest edge $\{a,b\}$ in $\alpha \cup
\beta$\,, the number of vertices greater than $b$ in $[2n-l-1]$ is at least $n-l-1$\,. This implies that $b \leq (2n-l-1)-(n-l-1) = n$\,. So $I_1$ is well-defined. Obviously, $I_1$ changes the number of elements of $\beta$ by one, i.e., $I_1$ is sign-reversing and double applications of $I_1$ on $(\alpha,\beta)$ will clearly restore $(\alpha,\beta)\,$. This construction cannot be accomplished, if both $\alpha$ and $\beta$ are empty sets, i.e., $n=k=l\,$. Thus we obtain the desired result. 
\end{proof}
\begin{eg}
Given $n=7$\,, $k=5$\,, $l=2$\,, let $\alpha=\{\,\{2,3\},\{4,7\}\,\}$ and $\beta=\{\,\{1,10\},\{5,11\},\{8,9\}\,\}$\,. Then $\alpha$ and $\beta$ are matchings in $K_{11}$ satisfying $(\alpha,\beta) \in U\,$. Take the smallest edge $e=\{2,3\}$ in $\alpha\cup\beta\,$. Since $e$ is an element of  $\alpha\,$, the image of $(\alpha,\beta)$ under $I_{1}$ is $(\alpha',\beta')\,$, where $\alpha'=\alpha\setminus\{e\}=\{\,\{4,7\}\,\}$ and $\beta'=\beta\cup\{e\}=\{\,\{1,10\},\{2,3\},\{5,11\},\{8,9\}\,\}\,$. It is easy to check that $I_{1}(\alpha',\beta')=(\alpha,\beta)$\,. (See Figure~\ref{fig:I_1}.)
\begin{figure}\label{fig:I_1}
\begin{center}
\includegraphics[scale=.7]{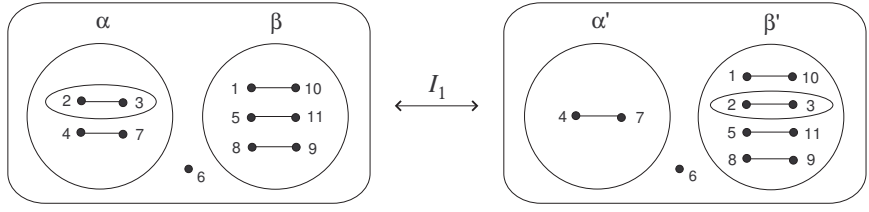}
\caption{The involution $I_1$\,.}
\end{center}
\end{figure}
\end{eg}

%%%%%%%%%%%%%%%%%%%%%%%%%%%%
\subsection{An involution for the second inverse formula}
\begin{thm}\label{thm:2ortho}
For all nonnegative integers $n$ and $l$\,, we have
$$\sum_{k=0}^{n}\,b(n,k)\,B(k,l) = \delta_{n,l} \,.$$
\end{thm}
\begin{proof}
Consider the complete graph $K_{2n-k}$\,. Let $V_k$ be the set of ordered pairs $(\alpha, \beta)$ such that
\begin{enumerate}
\item $\alpha$ is an $(n-k)$-matching in $K_{2n-k}$ in which  the vertex $2n-k$ is unsaturated under $\alpha$\,; and
\item $\beta$ is a $(k-l)$-matching in $K_{2n-k}$ such that $\alpha$ and $\beta$ saturate no vertex in common.
\end{enumerate}
Set $V=\bigcup V_k$\,. Define the sign of an element $(\alpha,\beta) \in V$ by ${\rm sgn}(\alpha,\beta) = (-1)^{|\alpha|}$\,. Then we get
\begin{equation*}
\sum_{(\alpha,\beta)\in V}{\rm sgn}(\alpha,\beta)
=\sum_{k}\,(-1)^{n-k}\,m(2n-k-1,n-k)\,\,m(k,k-l) = \sum_{k}\,b(n,k)\,B(k,l) \,.
\end{equation*}
To prove the theorem, it suffices to find a sign-reversing involution $I_2$ on $V$\,, which has exactly one fixed point with sign $+1$.

Suppose that $(\alpha, \beta) \in V_k$\,. Now consider the largest edge $e$ of $\alpha \cup \beta$ according to the colex order in $\alpha \cup \beta$\,. Define the map $I_2$ as follows:
\begin{itemize}
\item If $e \in \alpha$\,, then set $I_{2}(\alpha,\beta)=(\,\alpha\setminus\{e\}\,,\,\beta\cup\{e\}\,)$\,. Since $e \in \alpha$\,, neither $\alpha$ nor $\beta$ saturates the vertex $2n-k$\,. Moreover, the matching $\alpha\setminus\{e\}$ cannot saturate the vertex $2n-k-1$\,. So $(\,\alpha\setminus\{e\}\,,\,\beta\cup\{e\}\,)$ belongs to $V_{k+1}\,$.
\item If $e \in \beta$\,, then set $I_{2}(\alpha,\beta)=(\,\alpha\cup\{e\}\,,\,\beta\setminus\{e\}\,)$\,. Regard $(\,\alpha\cup\{e\}\,,\,\beta\setminus\{e\}\,)$ as a pair of matchings in $K_{2n-k+1}$\,. Clearly, the matching $\alpha\cup\{e\}$ cannot saturate the vertex $2n-k+1$\,. Thus $(\,\alpha\cup\{e\}\,,\,\beta\setminus\{e\}\,)$ belongs to $V_{k-1}\,$.
\end{itemize}
Evidently $(I_2 \circ I_2)(\alpha,\beta)=(\alpha,\beta)$ and $I_2$ is a sign-reversing map, whenever $\alpha \cup \beta \ne \emptyset$\,. The unique fixed point of $I_2$ is $(\alpha, \beta)\,$, where $\alpha$ and $\beta$ are both empty, i.e., $n=k=l\,$. This completes the proof.
\end{proof}
\begin{eg}
Given $n=10$\,, $k=8$\,, $l=5$\,, let $\alpha=\{\,\{2,3\},\{4,11\}\,\}$ and $\beta=\{\,\{1,7\},\{5,10\},\{8,9\}\,\}$\,. Then $\alpha$ and $\beta$ are matchings in $K_{12}$ satisfying $(\alpha,\beta) \in V_8\,$. Take the largest edge $e=\{4,11\}$ in $\alpha\cup\beta\,$. Since $e$ is an element of  $\alpha\,$, the image of $(\alpha,\beta)$ under $I_{2}$ is $(\alpha',\beta')\in V_9\,$, where $\alpha'=\alpha\setminus\{e\}=\{\,\{2,3\}\,\}$ and $\beta'=\beta\cup\{e\}=\{\,\{1,7\},\{4,11\},\{5,10\},\{8,9\}\,\}\,$. It is easy to check that $I_{2}(\alpha',\beta')=(\alpha,\beta)$\,. (See Figure~\ref{fig:I_2}.)
\begin{figure}\label{fig:I_2}
\begin{center}
\includegraphics[scale=.7]{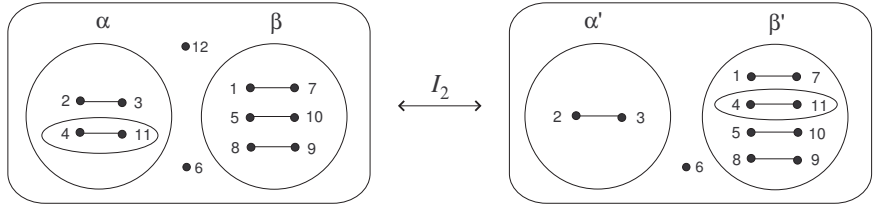}
\caption{The involution $I_2$\,.}
\end{center}
\end{figure}
\end{eg}

%%%%%%%%%%%%%%%%%%%%%%%
%%%%%%%%%%%%%%%%%%%%%%%
\section{Log-concavity of the Bessel numbers}
In this section, we show that  both Bessel numbers of the second kind $B(n,k)$ and signless Bessel numbers of the first kind $a(n,k)$ form log-concave sequences. In particular, we will check each result by constructing an explicit injection.

%%%%%%%%%%%%%%%%%%%%%%%%%%%%%%%%%%
\subsection{Log-concavity of the Bessel numbers of the second kind}  \label{subsec:log2}
\begin{thm} \label{thm:log2nd}
The Bessel numbers of the second kind $\{\,B(n,k)\,\}_{k \geq 0}$ form a log-concave sequence:
$$B(n,k-1)\cdot B(n,k+1)\leq B(n,k)^2 \quad \mbox{for all \,$k \geq 1$}\,.$$
\end{thm}
\begin{proof}
In fact, this theorem is a special case of the log-concavity of the sequence of matching numbers \cite{G, K}. For reader's convenience, we sketch here Krathenttaler's idea \cite{K}. Since $B(n,k)=m(n,n-k)$\,, it suffices to show that
\begin{equation*}
m(n,n-k+1)\cdot m(n,n-k-1)\leq m(n,n-k)^2 \,.
\end{equation*}
We will construct an injection $I_K$ from the set of all pairs $(\alpha_1,\alpha_2)$ of an $(n-k+1)$-matching $\alpha_1$ in $K_n$ and an $(n-k-1)$-matching $\alpha_2$ in $K_n$ into the set of all pairs $(\beta_1 ,\beta_2)$ of an $(n-k)$-matching $\beta_1$ in $K_n$ and an $(n-k)$-matching $\beta_2$ in $K_n \,$. Let an $(n-k+1)$-matching $\alpha_1$ in $K_n$ and an $(n-k-1)$-matching $\alpha_2$ in $K_n$ be given. Color the edges of $\alpha_1$ and $\alpha_2$ in blue and red, respectively.

Let $G=G[\alpha_1 \cup \alpha_2]$ denote the subgraph of $K_n$ induced by all edges of $\alpha_1$ and $\alpha_2$\,. Consider the connected components of $G$\,. Since $\alpha_1$ and $\alpha_2$ are matchings, each connected component of $G$ is either a cycle or a path. Note that all edges in the components of $G$ are colored alternately. So all cycles and paths of even length have an equal number of blue and red edges. Meanwhile, in case a path of odd length, the number of blue edges in the path differs from the number of red edges by one. We call a path of odd length a blue path, if it has more blue edges than red edges, and a red path, otherwise.

Let $R$ be the set of all red paths of $G$ and let  $r$ denote the cardinality of $R$\,. Since $G$ has two more blue edges than red, the number of all blue paths of $G$ should be $r+2$\,. Consider the set $C$ of all blue and red paths of $G$\,. Label the paths of $C$ from $1$ to $2r+2$ in the increasing order of the largest vertices in the paths.

Let $A=\{a_{1},\ldots,a_{r} \} \subset [2r+2]$\,, where $a_1 < \cdots < a_r$\,. Define $\phi(A)=A \cup\{a_{j}+1\}$\,, where $j$ is the largest $i$ for which $a_{i}-2i$ is minimal (assume $a_0=0\,$). Then the map $\phi$ is an injection from all $r$-element subsets of $[2r+2]$ into $(r+1)$-element subsets of $[2r+2]$\,. (See \cite[pp.~33--39~and~p.~54]{SW}.) 

Thus given the $r$-element subset $S$ of $[2r+2]$ corresponding to $R$, we have an $(r+1)$-element subset $\phi(S)=S\cup\{t\}$. Note that the element $t$ corresponds to a blue path, say $P$. If we exchange the colors in $P$, then $P$ becomes a red path, say $P'$. Now we consider the graph $G'=(G\setminus P) \cup P'$\,. Obviously, $G'$ has $n-k$ blue edges and $n-k$ red edges. 

Finally set  $\beta_1$ be the blue edges in $G'$ and $\beta_2$ be the red edges in $G'$. Let $I_{K}(\alpha_1,\alpha_2)=(\beta_1 ,\beta_2)\,$. It is easy to check that this defines the desired injection.
\end{proof}
\begin{eg}
Let $\alpha_1=\{\,\{1,2\}, \{3,4\}\, \{6,11\}\, \{7,12\}, \{13,14\}, \{16,17\}, \{19,20\}, \{21,22\}, \{23,24\}\}$ and $\alpha_2=\{\,\{2,3\},\{6,7\},\{8,9\},\{11,12\},\{14,15\},\{16,21\},\{19,24\}\,\}$\,. Then $\alpha_1$ and $\alpha_2$ are matchings in $K_{25}$\,. Color the edges of $\alpha_1$ in blue $(-)$ and the edges of $\alpha_2$ in red $(\cdots)$. There are three blue paths $1-2\cdots3-4$, $17-16\cdots21-22$ and $20-19\cdots24-23$\, and one red path $8\cdots9$\,, which are labeled by $1$, $3$, $4$ and $2$\,, respectively.  So we have $r=1$ and $S=\{2\}$\,. Since $\phi(\{2\})=\{2,3\}$\,, the blue path $P=17-16\cdots21-22$ is changed to the red path $P'=17\cdots16-21\cdots22$\,. From the graph $G'=(G\setminus P) \cup P'$\,, we can extract two $8$-matchings in $K_{25}$\,: $\beta_1=\{\,\{1,2\},\{3,4\}\,\{6,11\}\,\{7,12\},\{13,14\},\{16,21\},\{19,20\},\{23,24\}\,\}$ and $\beta_2=\{\,\{2,3\},\{6,7\},\{8,9\},\{11,12\},\{14,15\},\{16,17\},\{19,24\},\{21,22\}\,\}$\,. (See Figure~\ref{fig:I_K}.)
\begin{figure} \label{fig:I_K}
\begin{center}
\includegraphics[scale=.6]{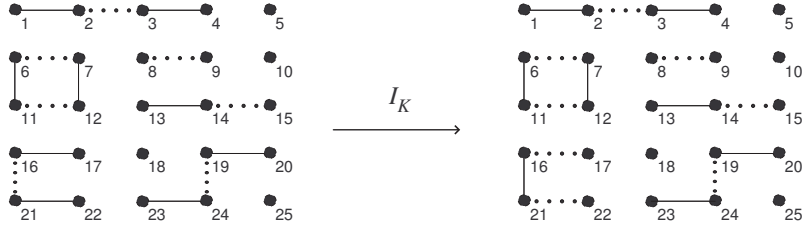}
\caption{The injection $I_{K}$.}
\end{center}
\end{figure}
\end{eg}

%%%%%%%%%%%%%%%%%%%%%%%%%%%%%%%%%%%%
\subsection{Log-concavity of the signless Bessel numbers of the first kind}
\begin{thm} \label{thm:log1st}
The sequence of signless Bessel numbers of the first kind $\{\,a(n,k)\,\}_{k \geq 0}$ is log-concave:
$$a(n,k-1) \cdot a(n,k+1) \leq a(n,k)^2 \quad \mbox{for all \,$k \geq 1$}\,.$$
\end{thm}

\begin{proof}
Since $a(n,k) = B(2n-k-1,n-1)=m(2n-k-1,n-k)$\,, this theorem is equivalent to the log-concavity of the sequence $m(2n-k-1,n-k)$\,, for $k=0, \ldots, n$\,. So it suffices to find an injection from the set of all pairs $(A_1,A_2)$ of an $(n-k+1)$-matching $A_1$ in $K_{2n-k}$ and an $(n-k-1)$-matching $A_2$ in $K_{2n-k-2}$ into the set of all pairs $(B_1,B_2)$ of an $(n-k)$-matching $B_1$ in $K_{2n-k-1}$ and an $(n-k)$-matching $B_2$ in $K_{2n-k-1}$\,.

Now we construct a new injection $I_{S}$ motivated by Krattenthaler's injection. Let an $(n-k+1)$-matching $A_1$ in $K_{2n-k}$ and an $(n-k-1)$-matching $A_2$ in $K_{2n-k-2}$ be given. Assign the colors blue to the edges of $A_1$ and red to the edges of $A_2$\,, respectively. Then there are two possible cases:
\begin{enumerate}
\item[(a)] The vertex $2n-k$ is unsaturated under $A_1$\,.
\item[(b)] The vertex $2n-k$ is saturated under $A_1$\,.
\end{enumerate}

In case (a), if we disregard the vertex $2n-k$ in $K_{2n-k}$ and regard $K_{2n-k-2}$ as the subgraph of $K_{2n-k-1}$ induced by the vertices in $[2n-k-2]$\,, then Krattenthaler's injection $I_K$ yields a blue $(n-k)$-matching $B_1$ in $K_{2n-k-1}$ and a red $(n-k)$-matching $B_{2}$ in $K_{2n-k-1}$\,. If the vertex $2n-k-1$ is saturated under $A_1$\,, then the blue path containing the vertex $2n-k-1$ should be labeled with the largest number $2r+2$\,, where $r$ is the number of red paths in $G[A_{1}\cup A_{2}]$\,. By the property of the map $\phi$ \cite[Chap.~2,~ex.~25]{SW}, the largest number $2r+2$ is not contained in $(\phi(S)\setminus S)$ for any $r$-subset $S$ of $[2r+2]$\,. So a blue path containing the vertex $2n-k-1$ is not selected to change into a red path. Moreover, the vertex $2n-k-1$ cannot be saturated under $A_2$\,. Therefore $B_2$ cannot saturate the vertex $2n-k-1$\,.

In case (b), we introduce a new mapping as follows. Consider the sets $X$ and $Y$ of all unsaturated vertices under $A_1$ and $A_2$\,, respectively. Let $x$ be the vertex which is adjacent to the vertex $2n-k$\,. Choose the element $y$ in $Y$ with the same relative order as $x$ in $X \cup \{x\}$\,. Since the cardinality of $X \cup \{x\}$ is $k-1$ and the cardinality of $Y$ is $k$, the selection $y$ is well-defined. Cut the blue edge $e_b =\{x,2n-k\}$ and join two vertices $y$ and $2n-k-1$ by the red edge $e_r$\,. Then $A_1\setminus \{e_{b}\}$ becomes a blue $(n-k)$-matching $B_1$ in $K_{2n-k-1}$ and $A_2 \cup \{e_{r}\}$ becomes a red $(n-k)$-matching $B_2$ in $K_{2n-k-1}$\,. Define $I_{N}(A_{1},A_{2})=(B_{1},B_{2})$\,. In this case, the vertex $2n-k-1$ is saturated under $B_2$\,. 

Let $I_{S}$ be the map $I_{K} \cup I_{N}$, i.e.,
\begin{equation*}
I_{S}(A_{1},A_{2})=\begin{cases}
\,I_{K}(A_{1},A_{2})\,, & \mbox{if $A_1$  does not saturate $2n-k$\,,}\\
\,I_{N}(A_{1},A_{2})\,,&\mbox{if $A_1$  saturates $2n-k$\,.}\end{cases}
\end{equation*}
Since both $I_{K}$ and $I_{N}$ are injective and the saturating condition of the vertex $2n-k-1$ makes  that the image of  $I_{K}$ and the image of $I_{N}$ are disjoint, the map $I_{S}$ is a desired injection. 
\end{proof}

\begin{eg}
Let $A_{1}=\{\,\{1,2\},\{3,4\},\{5,10\},\{6,11\},\{7,12\},\{13,14\},\{18,19\},\{20,25\},\{23,24\}\,\} \in \M(25,9)$ and $A_{2}=\{\,\{2,3\},\{6,7\},\{8,9\},\{11,12\},\{14,15\},\{17,18\},\{19,20\}\,\} \in \M(23,7)$\,. Color the edges of $A_1$ in blue $(-)$ and the edges of $A_2$ in red $(\cdots)$. Since the vertex  $25\,(=2n-k)$ is saturated under $A_1$\,, we should apply the injection $I_N$\,. In this case, $X=\{8,9,15,16,17,21,22\}$\,,  $x=20$ and $Y=\{1,4,5,10,13,16,21,22,23\}$\,. Since $20$ is the sixth smallest element in $X\cup\{x\}$\,, we should choose $16$ for $y$\,, which is the sixth smallest element in $Y$\,. Now cut the blue edge $\{20,25\}$ and join $16$ and $24$ with a red edge. From this, we can extract two $8$-matchings in $K_{24}$\,: $B_{1}=\{\,\{1,2\},\{3,4\},\{5,10\},\{6,11\},\{7,12\},\{13,14\},\{18,19\},\{23,24\}\,\}$ and $B_{2}=\{\,\{2,3\},\{6,7\},\{8,9\},\{11,12\},\{14,15\},\{17,18\},\{19,20\},\{16,24\}\,\}$\,. (See Figure~\ref{fig:I_N}.)
\begin{figure} \label{fig:I_N}
\begin{center}
\includegraphics[scale=.6]{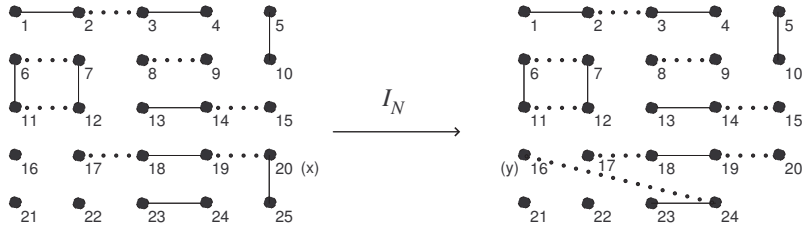}
\caption{The injection $I_{N}$\,. }
\end{center}
\end{figure}
\end{eg}
\begin{rmk}
Theorem~\ref{thm:log2nd} and~\ref{thm:log1st} imply the matching numbers  $\{m(n,k)\}$  are log-concave for both coordinates: 
\begin{eqnarray*}
m(n,k-1) \cdot m(n,k+1) \leq m(n,k)^2 \qquad &\mbox{for all}& k \geq 1\,, \\
m(n-1,k) \cdot m(n+1,k) \leq m(n,k)^2 \qquad &\mbox{for all}& n \geq 1\,. 
\end{eqnarray*}  
\end{rmk}

%%%%%%%%%%%%%%%
%%%%%%%%%%%%%%%
\section*{Acknowledgements}
The authors are grateful to their advisor Dongsu Kim for his advice and encouragement.
The second author thanks Ira M. Gessel for his helpful remarks and suggestions.
%%%%%%%%%%%%%
%%%%%%%%%%%%%

\end{document}